\documentclass{amsart}

\usepackage{amssymb}

\newtheorem{thm}{Theorem}

\newtheorem{lem}[thm]{Lemma}
\newtheorem{cor}[thm]{Corollary}

\newtheorem{prop}[thm]{Proposition}

\newtheorem{conj}[thm]{Conjecture}
   
\theoremstyle{definition}

\newtheorem{say}[thm]{}
\newtheorem{exmp}[thm]{Example}

\newtheorem{ques}[thm]{Question}    

\newtheorem{ack}{Acknowledgments}

\newtheorem{defn-thm}[thm]{Definition--Theorem}  
\newtheorem{defn-lem}[thm]{Definition--Lemma}  

\theoremstyle{remark}

\renewcommand{\c}[0]{{\mathbb C}}

\newcommand{\z}[0]{{\mathbb Z}}

\renewcommand{\r}[0]{{\mathbb R}} 

\renewcommand{\a}[0]{{\mathbb A}}

\newcommand{\p}[0]{{\mathbb P}}

\newcommand{\q}[0]{{\mathbb Q}}
\newcommand{\map}[0]{\dasharrow}
\newcommand{\qtq}[1]{\quad\mbox{#1}\quad}

\newcommand{\sing}[0]{\operatorname{Sing}}

\newcommand{\chow}[0]{\operatorname{Chow}}

\newcommand{\hilb}[0]{\operatorname{Hilb}}

\newcommand{\dio}[0]{\operatorname{dioph}}
\newcommand{\ddim}[0]{\operatorname{d-dim}}
\newcommand{\pole}[0]{\operatorname{pole}}
\newcommand{\Pole}[0]{\operatorname{Pole}}

\def\into{\DOTSB\lhook\joinrel\to}

\begin{document}
\bibliographystyle{amsalpha}

\title{Diophantine subsets of function fields of curves}
\author{J\'anos Koll\'ar}

\maketitle

\today

Let $R$ be a commutative ring. A subset $D\subset R$ is called
{\it diophantine} if there are polynomials 
$$
F_i(x,y_1,\dots,y_n)\in R[x,y_1,\dots,y_n]
$$
such that the system of equations
$$
F_i(r,y_1,\dots,y_n)=0\quad \forall i
$$
has a solution $( y_1,\dots,y_n)\in R^n$
iff $r\in D$.

Equivalently, if there is a (possibly reducible) algebraic variety
$X_R$ over $R$ and a morphism
$\pi:X_R\to  \a^1_R$ such that 
$D=\pi(X_R(R))$. In this situation we call
$$
\dio(X_R,\pi):=\pi(X_R(R))\subset R
$$
the {\it diophantine set} corresponding to $X_R$ and $\pi$.

A  characterization of diophantine subsets of
$\z$ was completed in connection with Hilbert's 10th problem,
 but a description of diophantine subsets of
$\q$ is still not known. In particular, it is not known
if $\z$ is a diophantine subset of $\q$ or not.
(See \cite{poo} or the volume \cite{h10} for surveys and many recent results.) 

In this paper we consider analogous questions
where $R=k(t)$ is a function field of one variable
and $k$ is an uncountable {\it large field} of
characteristic 0. That is, for any $k$-variety $Y$
with a smooth $k$-point, $Y(k)$ is Zariski dense. 
Examples of uncountable  large fields are
\begin{enumerate}
\item $\c$ or any uncountable algebraically closed field,
\item $\r$ or any uncountable real closed field,
\item $\q_p, \q((x))$ or the quotient field of
any   uncountable local Henselian domain.
\end{enumerate}

Roughly speaking, we show that for such fields,
a diophantine subset of $k(t)$ is either very small
or very large. The precise result is   somewhat technical,
but here are two easy to state consequences which served
as   motivating examples.

\begin{cor}\label{main.cor}
 Let $k$ be an uncountable  large field of
characteristic 0. Then $k[t]$ is not a
diophantine subset of $k(t)$.
\end{cor}

\begin{cor}\label{bog.cor}
 Let $k$ be an uncountable  large field of
characteristic 0 and $K_2\supset K_1\supset k(t)$  
finite field extensions. 
Then $K_1$ is  a
diophantine subset of $K_2$ iff $K_1=K_2$.
\end{cor}

The latter gives  a partial answer to a question of Bogomolov:
When is a subfield $K_1\subset K_2$ diophantine in $K_2$?

It is possible that both of these corollaries  hold for any field $k$.
Unfortunately, my method says nothing about countable fields.
The geometric parts of the proof (\ref{getZ.prop})  and (\ref{useZ.prop})
work for any 
uncountable   field, but the last step (\ref{pf.of.main}) uses in an essential
way that $k$ is large.
\medskip

We use two ways to measure how large a
diophantine set is.

\begin{say}[Diophantine dimension and polar sets]
Let $B$ be a smooth, projective, irreducible curve over $k$.
One can think of a rational function $f\in k(B)$
as a section  of the first projection $\pi_1:B\times \p^1\to B$.
 This establishes a
one-to-one correspondence 
$$
k(B)\cup\{\infty\}  \quad\leftrightarrow\quad 
\{\mbox{sections of $\pi_1:B\times \p^1\to B$}\}.
$$
Any section $\sigma:B\to B\times \p^1$ can be identified
with its image, which gives a point in the Chow variety
of curves of $B\times \p^1$. This gives an injection
$$
k(B)\cup\{\infty\}\quad\into \quad \chow_1(B\times \p^1).
$$

Let $U$ be a countable (disjoint) union of $k$-varieties and
$D\subset U(k)$ a subset. Define the
{\it diophantine dimension}   of $D$ over $k$ as the
smallest $n\in \{-1,0,1,\dots,\infty\}$
 such that $D$ is contained in a countable
union of irreducible  $k$-subvarieties of $U$ of dimension
$\leq n$. It is denoted by $\ddim_kD$.
Note that  $\ddim_kD=-1$ iff $D=\emptyset$ and
 $\ddim_kD\leq 0$ iff $D$ is countable.

In particular, we can talk about the diophantine dimension
of $\dio(X,f)\subset k(B)\subset \chow_1(B\times \p^1)$.

For $f\in k(B)$, let $\pole(f)$
denote its divisor of poles.
For $D\subset k(B)$ set 
$$
\Pole_n(D):=\{\pole(f): f\in D\mbox{ and } \deg\pole(f)=n\}.
$$
I think of $\Pole_n(D)$ as a subset of
the $\bar k$-points of the $n$th symmetric power $S^nB$.

Taking each point with multiplicity $r\geq 1$ gives
embeddings $S^mB\to S^{rm}B$, whose image I denote
by $r\cdot S^mB$.
\end{say}

With these definitions, the main result is the following
illustration of the ``very small
or very large'' dichotomy.

\begin{thm} \label{main.thm.alg} Let $k$ be an uncountable  large field of
characteristic 0 and $B$  a smooth, projective, irreducible curve over $k$.
Let  $X_{k(B)}$ be  a (possibly reducible) algebraic variety
of dimension $n$ 
 over $k(B)$ and
$\pi_{k(B)}:X_{k(B)}\to  \a^1_{k(B)}$  a  morphism.
Then
\begin{enumerate}
\item  either $\ddim_k \dio(X_{k(B)},\pi_{k(B)})\leq n$,
\item or $\ddim_k \dio(X_{k(B)},\pi_{k(B)})=\infty$ and
there  is a 0-cycle $P_a\in S^aB$  and $r>0$ such that
for every $m>0$ there is a smooth, irreducible
$k$-variety  $D_m$ and a morphism $\rho_m:D_m\to  S^{a+rm}B$ such that 
\begin{enumerate}
\item $D_m(k)\neq\emptyset$,
\item $\Pole_{a+rm}\bigl(\dio(X_{k(B)},\pi_{k(B)})\bigr)\supset 
\rho_m\bigl(D_m(k)\bigr)$,
 and
\item  the Zariski closure of
$\rho_m(D_m(k))$ contains
$P_a+r\cdot S^mB\subset S^{a+rm}B$.
\end{enumerate}
\end{enumerate}
\end{thm}

\begin{say}[Proof of the Corollaries]
In trying to write a subset $D\subset k(B)$ as $D=\dio(X_{k(B)},\pi_{k(B)})$,
we do not have an a priori bound on $\dim X_{k(B)}$, thus
the assertion $\ddim_k \dio(X_{k(B)},\pi_{k(B)})=\infty$
 is hard to use.
The Corollaries \ref{main.cor} and \ref{bog.cor}  both
follow from the more precise results about
the distribution of poles.

If $B=\p^1$, then a rational function with at least 2 poles on
$\p^1$ is not a polynomial, thus Theorem \ref{main.thm.alg}
implies Corollary \ref{main.cor}.

Next consider Corollary \ref{bog.cor}.
Let $K_1=k(B_1)\subsetneq K_2=k(B_2)$ be a degree $d>1$
 extension
of function fields of 
smooth, projective, irreducible $k$-curves.
By Riemann-Roch, any zero cycle of degree $\geq 2g(B_1)$ defined over
$k$ is the polar set of some $f\in k(B_1)$.
Pulling back gives  a map $j:S^mB_1\to S^{md}B_2$,
thus 
$$
\Pole_{n}(K_1) = 
\begin{cases}
 j\bigl((S^mB_1)(k)\bigr)\qtq{if $n=md\geq 2dg(B_1)$, and}\\
\emptyset  \qtq{if  $d\not\vert n$.}
\end{cases}
$$
If $b_1\neq b_2\in B_2$ map to the same point of $B_1$,
then a 0-cycle in $j\bigl(S^mB_1\bigr)$
contains either both $b_1$ and $b_2$ or neither.
Thus the Zariski closed set $j\bigl(S^mB_1\bigr)$
never contains a set of the form
$P_a+r\cdot S^mB_2$.
By (\ref{main.thm.alg}.2.c), this shows that
$K_1$ is not diophantine in $K_2$, proving
Corollary \ref{bog.cor}.
\end{say}

\begin{exmp} (1) 
The bound $n$ in (\ref{main.thm.alg}.1) is actually sharp,
as shown by the following.

Note first that any $k(t)$-solution of
$x^3+y^3=1$ is constant. Set
$$
X_n:=(x_1^3+y_1^3=\cdots=x_n^3+y_n^3=1)\subset 
\a^{2n}
$$
and 
$$
\pi:(x_1, y_1, \dots, x_n,y_n)\mapsto x_1+x_2t+\cdots+x_nt^{n-1}.
$$
Then 
$\dim X=n$ and for $k=\c$ or $k=\r$, 
$\dio(X_n,\pi)$ is the set of all degree $\leq n-1$ polynomials.

Using similar constructions
 one can see that any   (finite dimensional) $k$-algebraic subset
of $k(t)$ is diophantine when
$k$ is algebraically closed or real closed.
These are the ``small'' diophantine subsets of $k(t)$.

(2) The somewhat unusual looking condition
about the Zariski closure of $D_m$ in
(\ref{main.thm.alg}.2.c) is also close to being optimal.
For $g\in k(t)$  and $r>0$ consider the diophantine set
$$
L_{g,r}:=\{f\in k(t):\exists h\mbox{ such that } f= gh^r\}.
$$
Then, up to some lower dimensional contribution coming from
possible cancellations between poles and zeros of $g$ and $h^r$, 
$\Pole_{n}(L_{g,r})$  
equals
 $\pole(g)+r\cdot \bigl(S^mB\bigr)(k)$ if $n=\deg\pole(g)+rm$ and
$\emptyset$ otherwise.

\end{exmp}

\begin{say}
If $k=\c$ then our proof  shows
that in case (\ref{main.thm.alg}.2) 
there is a finite set $P\subset B(\c)$ such that
for every $p\in B(\c)\setminus P$ 
there is an  $f_p\in \dio(X,\pi)$
with a pole at $p$.

If $k=\r$, then we guarantee many poles, but one may not get any
real poles. 
To get examples, note that $h\in \r(t)$
is everywhere nonnegative on $\r$ iff $h$ is a sum of 2 squares.
Thus for any $g\in \r(t)$, the set
$$
L_1(g):=\{f\in \r(t): f(t)\leq g(t) \ \forall t\in \r\}
$$
is diophantine. $L(g)$ is infinite dimensional but if
$g\in \r[t]$ then no element of $L(g)$ has a real pole.

From the point of view of our proof   a more interesting
example is the diophantine set
$$
L_2(g):=\{f\in \r(t): \exists c\in \r,\ 
 f^2(t)\leq c^2\cdot g^2(t) \ \forall t\in \r\}.
$$
The elements of $L_2(g)$ are unbounded everywhere
yet 
no element of $L_2(g)$ has a pole in $\r$ if $g$ is a polynomial.

This leads to the following question.
\end{say}

\begin{ques}
Is $\r[t]_{\r}$, the set of all rational functions
without  poles   in $\r$, diophantine?
\end{ques}

There should be some even stronger variants
of the ``very small
or very large'' dichotomy,
 especially
over $\c$.  As a representative case, I propose the
following.

\begin{conj} Let $D\subset \c(t)$ be a diophantine subset
which contains a Zariski open subset of $\c[t]$.
(Meaning, for instance, that 
$D$ contains a Zariski open subset 
of the space of degree $\leq n$ polynomials for infinitely many $n$.)
Then $\c(t)\setminus D$ is finite.
\end{conj}

In connection with Bogomolov's question, I would hazard
the following:

\begin{conj} Let $k$ be a  large field and $K_1\subset K_2$ function fields
of $k$-varieties. Then $K_1$ is diophantine in $K_2$
iff $K_1$ is algebraically closed in $K_2$.
\end{conj}

\begin{say}
The proof  of (\ref{main.thm.alg})
relies on the theory of rational curves on algebraic varieties.
A standard reference is \cite{rc-book}, but non-experts may prefer
the more introductory lectures of \cite{ar-ko}. 

The proof is divided into three steps.

First we show that  if $\ddim_k \dio(X_{k(B)},\pi)\geq n+1$
then there is a rationally connected  
(cf.\ (\ref{rc.say})) subvariety $Z_{k(B)}\subset X_{k(B)}$
such that $\pi|_Z$ is nonconstant and $Z_{k(B)}$ has a smooth $k(B)$-point.
This relies on the {\it bend and break} method of
Mori \cite{mori1}. In a similar  context it was first used in
\cite{ghms}.

Then we show, using the {\it deformation of combs} 
technique developed in \cite{kmm2, rc-book, ghs, 0-spec},
 that for any such $Z_{k(B)}$, 
there are infinitely many $k$-varieties $S_m$ and maps
$S_m\times B\map Z_{k(B)}$
which give injections $S_m(k)\into Z_{k(B)}(k(B))$.

Both of these steps are geometric, but the
statements are formulated to work  over an arbitrary field $L$.

Finally, if $k$ is a large field, then 
each $S_m(k)$ is ``large'', which shows that
$Z_{k(B)}(k(B))$  is ``very large''.

For all three  steps it is better to replace $\pi:X_{k(B)}\to \a^1_{k(B)}$ with
a morphism of  $k$-varieties
  $f:X\to B\times \p^1$.
\end{say}

\begin{prop}\label{getZ.prop}
 Let $L$ be any field and $B$ a smooth, projective, irreducible curve
over $L$. Let $f:X\to B\times \p^1$ be an $L$-variety
of dimension $n+1$ and consider  the corresponding
diophantine set $\dio(X_{L(B)},f)\subset L(B)$.
Then
\begin{enumerate}
\item  either $\ddim_L \dio(X_{L(B)},f)\leq n$,
\item or there is a subvariety $Z\subset X$  such that
\begin{enumerate}
\item  $Z\to B\times \p^1$ is dominant, 
\item the generic fiber of $Z\to B$ is rationally connected, and
\item there is a rational 
section $\sigma:B\map Z$ whose image is not contained in $\sing Z$.
\end{enumerate}
\end{enumerate}
\end{prop}

\begin{prop}\label{useZ.prop}
 Let $L$ be an infinite field and $B$ a smooth, projective, irreducible curve
over $L$.
 Let $f:Z\to B\times \p^1$ be a smooth, projective $L$-variety
such that
\begin{enumerate}
\item $Z\to B\times \p^1$ is dominant,
\item the generic fiber of $Z\to B$ is separably rationally connected, and
\item there is a 
section $\sigma:B\to Z$.
\end{enumerate}
Then, for some $r>0$ and for all $m>0$ in an arithmetic progression,  
there are
\begin{enumerate}\setcounter{enumi}{3}
\item a smooth, irreducible $L$-variety $S_m$ 
with an $L$-point, and
\item a dominant rational map $\sigma_m:S_m\times B\map Z$
which commutes with projection to $B$
\end{enumerate}
 such that
the Zariski closure of the image of
$f\circ \sigma_m:S_m\map \chow_1(B\times \p^1)$ contains
$$
[f\circ \sigma(B)]+r[\{b_1\}\times \p^1]+\cdots + r[\{b_m\}\times \p^1]
\qtq{for every} b_i\in B(\bar L).
$$
\end{prop}

\begin{say}[Spaces of sections]
 Let $L$ be any field, $B$ a smooth, projective, irreducible curve
over $L$ and $f:X\to B$ a projective morphism.
A section of $f$ (defined over some $L'\supset L$)
can be identified with the corresponding $L'$-point in
the Chow variety of 1-cycles  $\chow_1(X)$.
All sections $\Sigma(X/B)$ defined over $\bar L$ form 
an open set of $\chow_1(X)$.
Indeed, if $H$ is an ample line bundle on $B$ of degree $d$
then a 1-cycle $C$ is a section iff $C$ is irreducible
(an open condition) and $(C\cdot f^*H)=d$
(an open and closed condition).
This procedure realizes  $X_{k(B)}\bigl(k(B)\bigr)$ as the set of $k$-points of
a countable union of
algebraic $k$-varieties $\Sigma(X/B)=\cup_i \Sigma_i$.

The choice of the $\Sigma_i$ is not canonical.
Given $X\to B$, we get ``natural'' irreducible components,
but for fixed generic fiber $X_{k(B)}$, these components
depend on the choice of $X$.  Any representation gives, however,
the same constructible sets.
We usually make a further decomposition.
Since every variety is a finite set-theoretic union of locally closed
smooth subvarieties, 
we may choose the  $\Sigma_i$ such that each one is smooth and irreducible.

As an explicit example, consider $B=\p^1$.
Then $k(B)\cong k(t)$ and every $f\in k(t)$
can be uniquely written (up to scalars) as 
$$
f=\frac{a_0+a_1t+\cdots+a_nt^n}{b_0+b_1t+\cdots+b_nt^n},
$$
where the nominator and the denominator are relatively prime
and at least one of  $a_n$ or $b_n$ is nonzero.
For any $n$, all such $f$ 
form an open subset
$$
\Sigma_n\subset \p(a_0:a_1:\cdots:a_n:b_0:b_1:\cdots:b_n)\cong \p^{2n+1}.
$$
\end{say}

\begin{say}[Very dense subsets]\label{zvd}
 Let $U$ be an irreducible variety
over a field $L$. We say that a subset $D\subset U(\bar L)$ is
{\it Zariski very dense} if $D$ is not contained in a countable
union of  $L$-subvarieties $V_i\subsetneq U$.

It is easy to see that for any $D$, there are
countably many closed, irreducible $L$-subvarieties $W_i\subset U$
such that 
$D\subset \cup_iW_i(\bar L)$ and 
$D\cap W_i(\bar L)$ is Zariski very dense in $W_i$ for every $i$.
There is a unique irredundant choice of these $W_i$.
\end{say}

\begin{say}[Proof of (\ref{getZ.prop})]
Write $X=\cup X_i$ as  a finite set-theoretic union of locally closed,
smooth, connected  varieties.
If (\ref{getZ.prop}) holds for each $X_i$ then it
also holds for $X$,
 thus
we may assume that $X$ is smooth and irreducible.
Let $X'\supset X$ be a smooth compactification such that
$f$ extends to $ f': X'\to B\times \p^1$.

As before,  there are countably many disjoint, irreducible, smooth 
$L$-varieties  $\cup_i\Sigma_i= \Sigma(X'/B)$ and morphisms
$u_i:B\times \Sigma_i\to X'$ commuting with projection to $B$
giving all $\bar L$-sections of $f'$.
As in (\ref{zvd}), there are countably many disjoint, irreducible, smooth 
$L$-varieties  $S_i\subset \Sigma(X'/B)$
such that each 
$S_i(L)$ is Zariski very dense in $S_i$ and
the $L$-sections of $X'\to B$
are exactly given by $\cup_i S_i(L)$.

Composing  $u_i$ with $f'$, we obtain maps
$$
f'_*:S_i\to \Sigma\subset\chow_1(B\times \p^1).
$$
There are 2 distinct possibilities. Either
\begin{enumerate}
\item  $\dim_L  f'_* (S_i)\leq n$ whenever
 $u_i(B\times S_i)\cap X\neq \emptyset$, or
\item  there is an $i_0$ such that 
$\dim_L  f'_* (S_{i_0})\geq n+1$
and $u_{i_0}(B\times S_{i_0})\cap X\neq \emptyset$.
\end{enumerate}
In the first case $\dio(X',f')$ is contained in
the union of the  constructible sets $f'_*(S_i)$, thus we  have
(\ref{getZ.prop}.1). This is always the case if $L$ is countable.

In the second case we construct $Z$ as required by
(\ref{getZ.prop}.2) using only the existence
of $u_{i_0}: B\times S_{i_0}\to X$.
Set $S:=S_{i_0}$ and $u:=u_{i_0}$. We can replace
$X'$ by a desingularization of
the closure of the image $u(B\times S)$. By shrinking $S$ we may assume
that $u$ lifts to $u:B\times S\to X'$.

For a  point $x\in X'$ let
$S_x\subset S$ be the subvariety
parametrizing those sections that pass through $x$.

Let us now fix $b\in B(\bar L)$ 
such that $u(\{b\}\times S)$ is dense in $X'_b$
and let $x$ run through
$X'_b$, the fiber of $X'$ over $b$. 
Since every section intersects $X'_b$,
$S=\cup_{x\in X'_b}S_x$ and so 
$f'_*(S)=\cup_{x\in X'_b}f'_*(S_x)$.
By assumption $\dim_L f'_*(S)\geq n+1$ and
$\dim_L X'_b=n$,
hence $\dim_L f'_*(S_x)\geq 1$
 for general $x\in X'_b(\bar L)$.
In particular, there is a 1-parameter family of sections
$C_x\subset S_x$ such that  
$$
f'\circ u :B\times C\to X\to B\times \p^1
$$
is a nonconstant family of sections passing
through the point $f'(x)$.

By (\ref{b-and-b}), this leads
to a  limit 1--cycle  of the form
$$
A+\{b\}\times \p^1+ (\mbox{other fibers of $\pi_1$})
$$
where $A$ is a section of $\pi_1:B\times \p^1\to B$.

Correspondingly, 
we get a limit 1--cycle in $X'$ of the form
$$
A_x+R_x+ (\mbox{other rational curves})
$$
where $A_x$ is a section of $X'\to B$ which 
 dominates $A$ and $R_x$ is a connected union of 
rational curves which  dominates $\{b\}\times \p^1$.
Note also that  $x\in R_x$.

Thus we conclude that for general $x\in X'_b(\bar L)$, there is
a connected union of rational curves $x\in R_x\subset X'_b$
which dominates $\{b\}\times \p^1$.

As in (\ref{mrc.say}), let us 
 take the relative MRC-fibration $f':X'\stackrel{w}{\map} W'\map B$.

For very general  $x\in X'(\bar L)$ let
$X'_x$ be the fiber of $w$ containing $x$.
By (\ref{mrc.say}), $X'_x$ is closed in $X'$
and every rational curve in $X'$ that intersects $X'_x$
is contained in $X'_x$. In particular, $R_x\subset X'_x$
and hence $X'_x$ dominates  $\{b\}\times \p^1$.

Let now $p\in S(L)$ be a general point and $C\subset X'$ the corresponding
section. 
By assumption $S(L)$ is Zariski dense in $S$, hence
we may assume that $w$ is smooth at the generic point of $C$.
Let $Z'\subset w^{-1}(w(C))$ 
be the unique irreducible component that dominates $C$
and $Z=Z'\cap X$. It satisfies all the required properties.\qed
\end{say}

\begin{say}[Bend-and-break for sections] \label{b-and-b}
(Cf.\ \cite{mori1}, \cite[Sec.II.5]{rc-book}, \cite[Lem.3.2]{ghms})

Let $h:Y\to B$ be a proper morphism onto a smooth projective curve $B$.
Let $C$ be a smooth curve and
$u :B\times C\to Y$
 a nonconstant family of sections passing
through a fixed point $y\in Y$.

Then  $C$ can not be a proper curve
and for a suitable point  $c\in \bar C\setminus C$
the corresponding  limit 1-cycle is of the form
$$
\Sigma_y=A_y+R_y,
$$
where $A_y$ is a section of $h$  (which need not pass through $y$)
and $R_y$ is a nonempty union of rational curves contained in finitely many
fibers of $h$. Furthermore, $A_y+R_y$ is connected
and $y\in R_y$.

This holds whether we take the limit in the Chow variety of 1-cycles,
in the Hilbert scheme or in the space of stable maps.

\end{say}

\begin{say}[Rationally connected varieties]\label{rc.say} 
(Cf.\ \cite{kmm2}, \cite[Chap.IV]{rc-book}, \cite[Sec.7]{ar-ko})

Let $k$ be a field
and $K\supset k$  an uncountable algebraically closed field.
A smooth projective $k$-variety $X$ is called
{\it rationally connected} or {\it RC}
if for every point pair $x_1,x_2\in X(K)$
there is a $K$-morphism $f:\p^1\to X$ such that
$f(0)=x_1$ and $f(\infty)=x_2$.
$X$ is called {\it separably rationally connected} or {\it SRC}
 if for every point  $x\in X(K)$
there is a $K$-morphism $f:\p^1\to X$ such that
$f(0)=x$ and $f^*T_X$ is an ample vector bundle.
(That is, a sum of positive degree line bundles.)
Furthermore,  $f:\p^1\to X$ can be taken to be an
 embedding if $\dim X\geq 3$.
It is known that SRC implies RC and the two notions are equivalent in
  characteristic 0.

We may not have any rational curves over $k$, but
we can work with the universal family of
these maps $f:\p^1\to X$.
Thus, if $\dim X\geq 3$ and $p\in X$ is a $k$-point,
then  there is an irreducible, smooth  $k$-variety
$U$  and a $k$-morphism
$G:U\times \p^1\to X$
such that 
\begin{enumerate}
\item $G(U\times \{0\})=p$,
\item  $G_u:\{u\}\times \p^1\to X_{\bar k}$
is an embedding for every $u\in U(\bar k)$, and
\item  $G_u^*T_X$ is ample for every $u\in U(\bar k)$.
\end{enumerate}

By \cite[Thm.1.4]{rc-loc}, if $k$ is large
then we can choose $U$ such that $U(k)\neq \emptyset$.
\end{say}

\begin{say}[MRC fibrations]\label{mrc.say}
(Cf. \cite{kmm2}, \cite[Sec.IV.5]{rc-book})

Let $K\supset k$ be as above.
Let $X$ be a  smooth projective $k$-variety and
$g:X\to S$ a $k$-morphism. There is a unique (up to birational maps)
 factorization
$$
g:X\stackrel{w}{\map} W  \stackrel{h}{\map} S
$$
such that 
\begin{enumerate}
\item  for general $p\in W(K)$, the fiber $w^{-1}(p)$ is closed in $X$ and
rationally connected, and
\item  for very general $p\in W(K)$
(that is, for $p$ in a countable intersection of
dense open subsets) 
every rational curve in $X(K)$ which intersects $w^{-1}(p)$
and maps to a point in $S$ 
is contained in $w^{-1}(p)$.
\end{enumerate}

The map $w:X\map W$ is called the (relative)
{\it maximal rationally connected fibration} or 
{\it MRC fibration} of $X\to S$.
Note that if $X$ contains very few rational curves
(for example, if $X$ is an Abelian variety or a K3 surface)
then $X=W$.
\end{say}

 \begin{say}[Proof of (\ref{useZ.prop})]

Here we essentially reverse the procedure of the first part.
Instead of degenerating a 1-parameter family of sections
to get a 1-cycle consisting of a section + rational curves,
we start with a section, add to it suitably chosen rational curves
and prove that this 1-cycle can be written as the limit of
sections in many different ways.

We assume that $Z$ is smooth, projective.
If necessary, we take its product with $\p^3$
to achieve that $\dim Z\geq 4$.
This changes the space of sections $\Sigma(Z/B)$
but it does not change the image of $\Sigma(Z/B)$
in $L(B)$.

Apply (\ref{rc.say}) to $X= Z_{L(B)}$ and the point
$p= \sigma(B)$ to get
$$
G:U_{L(B)}\times \p^1\to Z_{L(B)}.
$$
Next replace $U_{L(B)}$ by  an $L$-variety $\tau:U\to B$
such that $G$ extends to $g:U\times \p^1\to Z$. 
By shrinking $U$ if necessary, we may assume that
for general $b\in B(\bar L)$, the corresponding
$$
g_b:U_b\times \p^1\to Z_b
$$
is a family of smooth rational curves passing through $\sigma(b)$
and $g_{b,u}^*T_{Z_b}$ is ample for every $u\in U_b$
where $g_{b,u}$ is the restriction of $g_b$ to $\{u\}\times \p^1$.

Given  distinct points $b_1,\dots, b_m\in  B(\bar L)$, 
let $B(b_1,\dots, b_m)$ be the comb assembled from
$B$ and $m$ copies of $\p^1$ where we attach
$\p^1_i$ to $B$ at $b_i$ (cf.\ (\ref{comb.say})).

By \cite[Thm.16]{0-spec}, there are 
$b_1,\dots, b_{m_0}\in  B(\bar L)$ and an embedding
$$
\sigma(g_{1},\dots, g_{m_0}):B(b_1,\dots, b_{m_0})\to Z
$$
given by $\sigma$ on $B$ and by $g_i:=g_{b_i,u_i}$
on the $\p^1_i$ for some $u_i\in U_{b_i}$ such that  the image, denoted by
$B(g_{1},\dots, g_{m_0})\subset Z$
is defined over $L$ and 
its normal bundle is as positive as one wants.
In particular, by (\ref{defcomb.say}), $B(g_{1},\dots, g_{m_0})$
gives a smooth point of the Hilbert scheme of $Z$.
Furthermore, for any further distinct  points
$b_{m_0+1},\dots, b_{m}$ and $g_{i}$ for $i=m_0+1,\dots,m$,
the resulting 
$$
\sigma(g_{1},\dots, g_{m}):B(b_1,\dots, b_{m})\to Z
$$
also gives a smooth point of the Hilbert scheme of $Z$.

Let $S_m$ denote the smooth locus of  the corresponding 
$L$-irreducible component
of the Hilbert scheme of $Z$. 
 $B(g_{1},\dots, g_{m})$ gives an $L$-point of $S_m$,
hence $S_m$  is geometrically irreducible.
By (\ref{defcomb.say}) 
the general point of $S_m$ corresponds to a section
of $f$, 
the universal family is a product
over a dense open subset of $S_m$ and 
we have a dominant rational map
$\sigma_m:S_m\times B\map Z$.

For a given $m$, it is not always possible to choose
$b_{m_0+1},\dots, b_m$ such that the set
$b_{1},\dots, b_m$ is defined over $L$. To achieve this,
choose a generically finite and dominant map
$\rho:U\map \a^{\dim U}_L$. For general
$c\in \a^{\dim U}(L)$, its preimage
$\rho^{-1}(c)$ gives $\deg\rho$ general points in $U$
which are defined over $L$.  Thus we can choose
$b_{1},\dots, b_m$ to be defined over $L$
whenever $m-m_0$ is a multiple of $\deg\rho$.

Let us now consider 
$f_*(S_m)\subset   \chow_1(B\times \p^1)$.
By construction, it contains
$$
f_*\bigl(B(g_{1},\dots, g_{m})\bigr)=A+r\bigl(\{b_1\}\times \p^1\bigr)+
\cdots +r\bigl(\{b_m\}\times \p^1\bigr)
$$
where $A=f\circ\sigma(B)$ is a section of $B\times \p^1$ and  $r\geq 1$
is the common (geometric) degree of the maps 
$$
f\circ g_{b,u}:\{u\}\times \p^1\to \{b\}\times \p^1\subset B\times \p^1.
$$

The combs 
$B(g_{1},\dots, g_{m})\subset Z$ have some obvious deformations
where we keep $B$ fixed and vary the points $b_i$ and the maps $g_{i}$
with them.
By construction, the points $b_i$ can vary independently.
The images of these 1-cycles in $B\times \p^1$
are of the form
$$
A+r\bigl(\{b'_1\}\times \p^1\bigr)+
\cdots +r\bigl(\{b'_m\}\times \p^1\bigr),
$$
where the $b'_i$ vary independently.
\qed
\end{say}

\begin{say}[Combs]\label{comb.say}
A  {\it comb} assembled from
$B$ and $m$ copies of $\p^1$ attached
 at the points $b_i$  is a curve  $B(b_1,\dots, b_m)$  obtained from
$B$ and $\{b_1,\dots, b_m\}\times \p^1$
by identifying the points $b_i\in B$ and $(b_i,0)\in \{b_i\}\times \p^1$.

A comb can be pictured as below:
$$
\begin{array}{c}
\begin{picture}(100,100)(40,-70)
\put(0,0){\line(1,0){180}}
\put(20,10){\line(0,-1){60}}
\put(50,10){\line(0,-1){60}}
\put(160,10){\line(0,-1){60}}
\put(10,5){$b_1$}
\put(40,5){$b_2$}
\put(150,5){$b_m$}

\put(20,0){\circle*{3}}
\put(50,0){\circle*{3}}
\put(160,0){\circle*{3}}

\put(15,-60){$\p^1_1$}
\put(45,-60){$\p^1_2$}
\put(155,-60){$\p^1_m$}

\put(-25,-5){$B$}

\put(90,-30){$\cdots\cdots$}

\end{picture}\\
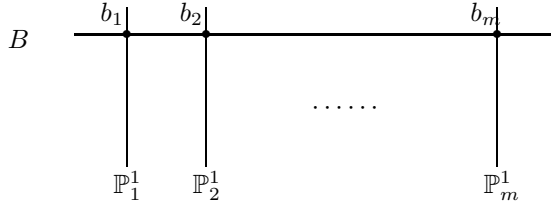
\mbox{Comb with $m$-teeth}
\end{array}
$$
\end{say}

\begin{say}[Deformation of reducible curves]\label{defcomb.say}
Let $X$ be a smooth projective variety and $C\subset X$
a connected curve with ideal sheaf $I_C$. Assume that $C$ has only nodes
as singularities.  By the smoothness criterion of the  Hilbert scheme 
\cite[p.221-21]{gro-fond}, if 
$H^1(C, \bigl(I_C/I_C^2\bigl)^*)=0$ then 
$[C]$ is a smooth point of the Hilbert scheme $\hilb(X)$
and there is a unique irreducible component $H_C\subset \hilb(X)$
containing $[C]$. Let $U_C\to H_C$ be the universal family.

If, in addition, $(I_C/I_C^2\bigl)^*$ is generated by global sections,
then a general point of $H_C$ corresponds to a smooth curve
and the natural map $U_C\to X$ is dominant.
\end{say}

\begin{say}[Proof of (\ref{main.thm.alg})]\label{pf.of.main}

Let us start with the $k(B)$-variety $X_{k(B)}$.
We can write it as the generic fiber of a quasi-projective $k$-variety
$X\to B$ and extend $\pi_{k(B)}$  to $f:X\to B\times \p^1$.
If (\ref{main.thm.alg}.1) fails then using (\ref{getZ.prop})
we obtain $Z\subset X$.
Take a compactification $\bar Z$ and a resolution
$Z_1\to \bar Z$ such that  the composite map
$Z_1\map B\times \p^1$ is a morphism. Next apply 
 (\ref{useZ.prop}) to $Z_1\to B\times \p^1$. We obtain,  for every $m$ 
in an arithmetic progression, 
 a dominant family of sections $\sigma_m:S_m\times B\to Z_1$.

There is a dense open subset $D_m\subset S_m$ such that
for every $q\in D_m(\bar k)$,
\begin{enumerate}
\item  the 
corresponding section   $\sigma_m(\{q\}\times B)\subset Z_1$ intersects $Z$,
and 
\item  the 
corresponding rational function $f\circ \sigma_m:\{q\}\times B\to \p^1$
has exactly $a+rm$ poles where
$a$ is the number of poles of $f\circ \sigma(B)$. 
\end{enumerate}
Thus the  the composite map
$$
\rho_m:=\pole\circ f\circ \sigma_m: D_m\to S^{a+rm}B
$$
 is defined. The condition
(\ref{main.thm.alg}.2.b) holds by construction
and the Zariski closure of
$\rho_m(D_m)$ contains
$P_a+r\cdot S^mB\subset S^{a+rm}B$ by (\ref{useZ.prop}.5),
where $P_a$ denotes the polar divisor of the
section $\sigma$, that is, the $0$-cycle
$f\circ \sigma(B)\cap \bigl(B\times \{\infty\}\bigr)$.

$S_m$ has a smooth $k$-point by (\ref{useZ.prop})
and $k$-points are Zariski dense since $k$ is a large field.
Thus $D_m(k)$ is Zariski dense in $D_m$. 
This implies both (\ref{main.thm.alg}.2.a) and (\ref{main.thm.alg}.2.c).

Finally, $D_m(k)$ is Zariski very dense in $D_m$
by (\ref{card.of.X}) and
$$
\ddim_k \dio(X,\pi)\geq \ddim_k \rho_m(D_m(k))= \dim_k \rho_m(D_m)
\geq \dim_k S^mB=m,
$$
where the middle  equality holds by (\ref{card.of.X}).
Thus $\ddim_k \dio(X,\pi)=\infty$. 

The only remaining issue is that our $m$ runs through an
arithmetic progression and is not arbitrary. 
If the progression is $m=b+m'c$ then
$a+r(b+m'c)=(a+rb)+(rc)m'$,
so by changing $a\mapsto a+rb, r\mapsto rc$
we get   (\ref{main.thm.alg}).
\qed

\end{say}

\begin{lem}\label{card.of.X} Let $X$ be a smooth and irreducible variety over
a large field $k$ such that $X(k)\neq\emptyset$.
Then $X(k)$ is not contained in the union of
fewer than $|k|$  subvarieties  $ X_{\lambda}\subsetneq X$. 
In particular, if $k$ is uncountable then 
 $X(k)$ is Zariski very dense in $X$.
\end{lem}

Proof.  Assume to the contrary that $X(k)\neq\emptyset$
but $X(k)\subset \cup_{\lambda\in \Lambda} X_{\lambda}$
where $|\Lambda|<|k|$ and $X_{\lambda}\neq X$.

If $\dim X\geq 2$, then pick 
 $p\in X(k)$ and let $\{H_t:t\in \p^1_k\}$ be a general pencil
of hypersurface sections of $X$ passing through $p$.
Since $|\Lambda|<|k|$, there is an $H_t$ such that
$H_t$ is smooth at $p$ and $H_t\not\subset X_{\lambda}$
for every $\lambda$. Thus
$H_t(k)\subset \cup_{\lambda\in \Lambda} (H_t\cap X_{\lambda})$
is a lower dimensional counter example.
Thus it is enough to prove (\ref{card.of.X})
when $X$ is a curve. Then lower dimensional $k$-subvarieties are just
points, thus we need to show that $|X(k)|=|k|$.

Let $\{H_t:t\in \p^1\}$ be a linear system on $X\times X$
which has $(p,p)$ as a base point and whose general member
is smooth at $(p,p)$. Since $k$ is large, each $H_t$ 
contains a $k$-point different from $(p,p)$.
Thus $|X(k)|=|X(k)\times X(k)|\geq |k|$.
\qed 

 \begin{ack}  I thank F.\ Bogomolov, 
K.\ Eisentr\"ager,  B.\ Poonen and J.\ Starr
for useful conversations and the referee for several insightful
suggestions. 
Partial financial support  was provided by  the NSF under grant number 
DMS-0500198. 
\end{ack}

\bibliography{refs}

\vskip1cm

\noindent Princeton University, Princeton NJ 08544-1000

\begin{verbatim}kollar@math.princeton.edu\end{verbatim}

\end{document}